\documentclass[aop,MSNbibl,seceqn,dvips]{arximspdf}
\usepackage{url,breakurl}
\doi{10.1214/13-AOP901} %kopijuoti is PTS
\volume{43}
\issue{1}
\pubyear{2015}
\firstpage{166}
\lastpage{187}

\makeatletter

\def\cal{\mathcal}

\def\eps{\varepsilon}
\def\P{{\mathbb P}}
\def\esp{{\mathbb E}}
\def\var{{\mathbb V}}
\def\Z{{\mathbb Z}}
\def\N{{\mathbb N}}
\def\R{{\mathbb R}}

\def\a{\alpha}
\def\l{\lambda}
\def\G{\Gamma}
\def\O{\Omega}

\def\D{{\cal D}}

\def\ra{\rightarrow}

\newtheorem{theo}{Theorem}[section]
\newtheorem{prop}[theo]{Proposition}
\newtheorem{lm}[theo]{Lemma}

\newproclaim{rmk}[theo]{Remark}
\newproclaim{df}[theo]{Definition}
\makeatother

\begin{document}
\begin{frontmatter}

\title{Differentiating the entropy of random walks on hyperbolic groups}
\runtitle{Differentiating the entropy of random walks}

\begin{aug}
\author[A]{\fnms{P.} \snm{Mathieu}\corref{}\ead[label=e1]{pierre.mathieu@cmi.univ-mrs.fr}}
\runauthor{P. Mathieu}
\affiliation{Aix Marseille Universit\'e}
\address[A]{Aix Marseille Universit\'e\\
CNRS, Centrale Marseille\\
LATP, UMR 7353\\
13453 Marseille\\
France\\
\printead{e1}} %adresu isvedimo komanda gale!
\end{aug}

% HISTORY:
\received{\smonth{9} \syear{2012}}
\revised{\smonth{11} \syear{2013}}

% ABSTRACT
%
\begin{abstract}
We show that the asymptotic entropy of a random walk on a nonelementary
hyperbolic group, with symmetric and bounded increments, is
differentiable and we identify its derivative as a correlation. We also
prove similar results for the rate of escape.
\end{abstract}

% KEYWORDS

\end{frontmatter}

%s1 #&#
\section{Introduction}

Let $\G$ be an infinite, countable, discrete group with neutral element
$\mathrm{id}$, and let $\mu$ be a probability measure on $\Gamma$.
The entropy of $\mu$ is defined as
\[
H(\mu):=\sum_{x\in\Gamma} \bigl(-\log\mu(x)\bigr) \mu(x).
\]
Note that $H(\mu)$ is nonnegative and may be infinite.

Let $\mu^n$ denote the $n$th convolution power of $\mu$.

We assume that $H(\mu)<\infty$.
It is easy to check that the sequence\break $(H(\mu^n))_{n\in\N}$ is
subadditive so that the following limit does exist:
%
%e1.1 #&#
\begin{equation}
\label{eq:0} h(\mu):=\lim_{n\ra\infty}\frac{1} n H\bigl(
\mu^n\bigr).
\end{equation}
The quantity $h(\mu)$ is called the \emph{asymptotic entropy} of
$\mu$.

The notion of asymptotic entropy was introduced by Avez in \cite
{kn:avez} in relation with random walk theory. In \cite{kn:avez2}, Avez
proved that, whenever $h(\mu)=0$, then $\mu$ satisfies the Liouville
property: bounded, $\mu$-harmonic functions are constant.
The converse was proved later; see \cite{kn:derrie} and \cite{kn:kaiver}.

Consider the random walk on $\Gamma$ whose increments are distributed
according to~$\mu$. The Liouville property is equivalent to the
triviality of the asymptotic $\sigma$-field of the random walk (its
so-called Poisson boundary); see \cite{kn:derrie} and \cite{kn:kaiver}
again. In more general terms, the entropy plays a central role in the
identification of the Poisson boundary of random walks in many
examples. We refer in particular to~\cite{kn:kaim} for groups with
hyperbolic features. In this latter case, the asymptotic entropy is
also related to the geometry of the harmonic measure through a
``dimension-rate of escape-entropy'' formula; see \cite{kn:bhm2} and the
references quoted therein.

The question of the regularity of $h$ as a function of $\mu$ was raised
by Erschler and Kaimanovich in \cite{kn:EK}, where it is proved
that, still for hyperbolic groups, the asymptotic entropy is
continuous. If we restrict ourselves to measures $\mu$ with fixed
finite support, Ledrappier recently proved in \cite{kn:led2} that
$h$ is Lipschitz continuous. We shall complement the result of
\cite{kn:led2} showing that, for a hyperbolic group $\Gamma$ and
restricting ourselves to symmetric measures
$\mu$ with a fixed finite support, the asymptotic entropy is
differentiable (Theorem~\ref{theo:entrop}).

There is an analogy between the definition of $h$ and the definition of
the \emph{rate of escape} of the random walk in some left-invariant
metric. More precisely, it can be proved that $h(\mu)$ coincides with
the rate of escape of the corresponding random walk in the so-called
Green metric; see the definition in Section \ref{sec:defre}. We also give
sufficient conditions on a metric ensuring that the rate of escape is
differentiable (Theorem~\ref{theo:rescp}).

Our approach completely differs from the one in \cite{kn:led2}. We
start with the simple observation that the derivative of the mean
position of a random walk is a correlation; see Section \ref{sec:heu}.
Thus, the natural candidate to be the derivative of the rate of escape
is some asymptotic covariance. These heuristics suggest a close
connection between the differentiability of the rate of escape and the
central limit theorem and explain the statement of Theorem~\ref{theo:rescp}.

As for the entropy, one deduces the differentiability of $h$ and the
value of its derivative from Theorem~\ref{theo:rescp} when choosing the
right metric (namely the Green metric) and observing that fluctuations
of the Green metric are of lower order---a fact that follows from the
``fundamental inequality'' between entropy and rate of escape and which
is true for random walks on nonamenable groups in general; see Section
\ref{sec:entrop}.

The version of the central limit theorem we need is a straightforward
extension of \cite{kn:bjork}. We also rely a lot on the ``Green metric
machinery'' developed in \cite{kn:bhm2}.

Let us emphasize that we do not claim that our results are optimal. It
is quite possible that the entropy and rate of escape are much more
regular that differentiable. It is actually known that the entropy and
rate of escape are analytic on the free group \cite{kn:led1} and that
the rate of escape is analytic for some Fuchsian groups~\cite{kn:HMM}.
One may hope analyticity holds for general hyperbolic groups (although
it does not hold for all groups, see \cite{kn:mm}). Anyway, we believe
the interpretation of the derivative as a correlation is rather
satisfactory, at least from an intuitive point of view. It clarifies
the connection between the regularity of the rate of escape and the
central limit theorem, an observation that seems to be new in our context.

Let us finally mention that the interpretation of the derivative of a
steady state (whatever it may mean) as some kind of correlation is a
well-known idea in theoretical physics or dynamical systems, where it
is sometimes called ``linear response theory'' or
``fluctuation-dissipation theory.'' See \cite{kn:ruelle} and other papers
of the same author.

%%%%%%%%%%%%%%%%%%%%%%%%%%%%%%%%%%%%%%%%%%%%%%%%%%%%%%%%%%%%%%%%%%%%%%%%%%%%%%%%%%%%%
%%%%%%%%%%%%%%%%%%%%%%%%%%%%% DEFINITIONS AND RESULTS
%%%%%%%%%%%%%%%%%%%%%%%%%%%%%%%%%%%%%%%%%%%%%%%
%%%%%%%%%%%%%%%%%%%%%%%%%%%%%%%%%%%%%%%%%%%%%%%%%%%%%%%%%%%%%%%%%%%%%%%%%%%%%%%%%%%%%

%s2 #&#
\section{Definitions and results}
\label{sec:defre}

%s2.1 #&#
\subsection{Definitions}

Let $\G$ be an infinite, countable, discrete group with neutral
element $\mathrm{id}$.
Let $d$ be a left-invariant proper metric on $\G$. We assume that $\G$
is finitely generated and that $d$ is equivalent to a word metric.

When the context is clear, we may also use the notation $\vert x-y\vert
$ to denote the distance
between $x$ and $y$ and $\vert x\vert=\vert x-\mathrm{id}\vert$.

We define the Gromov product of points $x,y\in\G$ with respect to the
base point $w\in\G$ by
\[
(x,y)_w:=\tfrac{1}2 \bigl(\vert x-w\vert+\vert y-w\vert-\vert x-y
\vert\bigr).
\]

We recall that the distance $d$ is called \emph{hyperbolic} if there
exists a constant $\tau$ such that
%
%e2.1 #&#
\begin{equation}
\label{eq:hyper} (x,y)_w\geq\min\bigl\{(x,z)_w,(z,y)_w
\bigr\} -\tau
\end{equation}
for all $x,y,z,w\in\Gamma$.

The group $\G$ is called hyperbolic if any (equivalently some) word
metric is hyperbolic.
A hyperbolic group is called \emph{nonelementary} if it is nonamenable
(which turns out to be equivalent to requiring $\G$ is not a finite
extension of~$\Z$).

We refer to \cite{kn:gh} for background material on hyperbolic groups.

From now on, we will assume that $\G$ is hyperbolic and nonelementary.

Following \cite{kn:bhm2}, we let $\D(\G)$ be the set of left-invariant
proper metrics on $\G$ which are both
equivalent to a word metric and hyperbolic.
Note that these last two conditions are not redundant. Indeed there
always exist nonhyperbolic (nongeodesic) metrics
equivalent to any word metric on $\Gamma$; see \cite{kn:bhm2},
Proposition~2.3. There may also exist hyperbolic metrics on $\G$ that
are not equivalent to a word metric.

We shall consider the following two compactifications of $\G$.
The \emph{visual \textup{(}Gromov\textup{)} compactification} is obtained by considering
all infinite sequences in $\G$ and identifying two such sequences,
say $(x_n)$ and $(y_n)$, if the Gromov product $(x_n,y_n)_w$ tends to
infinity. The \emph{horo-function \textup{(}Busemann\textup{)} compactification} is constructed
by identifying a point $x\in\G$ with the horo-function $k_x\dvtx \G
\rightarrow\R$, $k_x(y)=\vert y-x\vert-\vert w-x\vert$ and taking
the closure
for the topology of pointwise convergence. By Ascoli's theorem, this is
indeed a compact. The group $\G$ acts by homeomorphisms
on both compactifications. Up to equivariant homeomorphisms, the Gromov
compactification does not depend on the choice of
either the base point $w$ or the choice of $d\in\D(\G)$. The Busemann
compactification is also independent of the choice
of the base point but not of the choice of the metric. We shall say
that $d$ satisfies \emph{Assumption \textup{(}BA\textup{)}} if, up to equivariant
homeomorphisms, the Gromov and Busemann compactifications coincide.

Assumption (\emph{BA}) is in particular satisfied by the class of metrics
called \emph{Green metrics}. These are constructed as follows.
We call a probability measure $\mu$ on $\G$ ``symmetric'' if $\mu
(x^{-1})=\mu(x)$ for all $x\in\G$.
The support of $\mu$ is the set of $x\in\G$ such that $\mu(x)$ is
not zero.

Let $\mu$ be a probability measure on $\G$. We assume that $\mu$ is
symmetric and
that the support of $\mu$ is finite and generates the whole group $\G$.
The Green function associated to $\mu$ is defined by
\[
G^\mu(x)=\sum_{n=0}^\infty
\mu^n(x),
\]
where $\mu^n$ is the $n$th convolution power of $\mu$.

We assumed that $\G$ is nonamenable. Therefore, the sequence $\mu^n(x)$
exponentially converges to zero so that the series defining $G^\mu$
does converge.
The Green distance between points $x$ and $y$ in $\G$ is then
%
%e2.2 #&#
\begin{equation}
\label{eq:green}d_G^\mu(x,y):=\log G^\mu(\mathrm{id})-\log
G^\mu\bigl(x^{-1}y\bigr).
\end{equation}
In \cite{kn:bhm2}, we proved that $d_G^\mu$ belongs to $\D(\G)$.
Observe that $d_G^\mu$ need not be geodesic.
(As a matter of fact, it is not so difficult to deduce that $d_G^\mu$
is equivalent to a word metric from the nonamenability of $\G$.
That $d_G^\mu$ is hyperbolic is equivalent to a certain
multiplicativity property of the Green function $G^\mu$ which is the content
of Ancona's classical---and difficult---theorem; see \cite{kn:anc}
and the proof in \cite{kn:woess}.)

We now give the definition of the \emph{random walk} associated to a
probability measure $\mu$ on $\G$.
Because we will eventually use Radon--Nikodym transforms, it will be
more convenient to work with the canonical
construction on the set of trajectories on $\G$, say $\Omega=\G^\N$,
where $\N=\{0,1,\ldots\}$.
Given $\omega=(\omega_0,\omega_1,\ldots)\in\Omega$ and $n\in\N$,
we define the maps
$Z_n$ and $X_n$ from $\Omega$ to $\G$ by $Z_n(\omega):=\omega_n$,
and $X_n(\omega):=(Z_{n-1}(\omega))^{-1}Z_n(\omega)$. Thus,
$Z_n(\omega)$ gives the
position of the trajectory $\omega$ at time $n$, while
$X_n(\omega)$ gives its increment also at time $n$.
Following the usual usage in probability theory, we often omit to
indicate that random functions, as $Z_n$ or $X_n$, depend on $\omega$.

We equip $\Omega$ with the product $\sigma$-field (i.e., the smallest
$\sigma$-field for which all functions $Z_n$ are measurable).
The law of the random walk with increments distributed like $\mu$ is,
by definition, the unique probability measure on $\Omega$
under which $Z_0=\mathrm{id}$ and the random variables $(X_n)_{n\in\N}$ are
independent and distributed like $\mu$. We denote it with $\P^\mu$.
We also use the notation $\esp^\mu$ to denote the expectation with
respect to $\P^\mu$. Observe that the law of $Z_n$
under $\P^\mu$ is $\mu^n$.

We recall that, given a probability measure $\mu$ with finite support
and given a left-invariant metric $d\in\D(\G)$, the \emph{rate of escape}
of $\mu$ in the metric $d$ is defined by
\[
\ell(\mu;d):=\lim_{n\ra\infty}\frac{1} n \sum
_{x\in\G} d(\mathrm{id},x)\mu^n(x)
\]
and the \emph{asymptotic entropy} of $\mu$ is defined by
\[
h(\mu):=\lim_{n\ra\infty}\frac{1} n \sum
_{x\in\G} \bigl(-\log\mu ^n(x)\bigr) \mu
^n(x).
\]

Kingman's subadditive theorem implies that
\[
\ell(\mu;d)=\lim_{n\ra\infty}\frac{1} n d(\mathrm{id},Z_n)\quad
\mbox{and}\quad h(\mu )=\lim_{n\ra\infty}-\frac{1} n \log
\mu^n(Z_n),
\]
where both limits hold $\P^\mu$ almost surely as well as in $L^1(\O
,\P
^\mu)$.

%%%%%%%%%%%%%%%%%%%%%%%%%%%%%%%%%%%%%%%%%%%%%%%%%%%%%%%%%%%%%%%%%%%
%s2.2 #&#
\subsection{Differentiability of \texorpdfstring{$\ell$}{$ell$} and $h$}

In the sequel, we shall study the derivatives of the rate of escape and
the entropy of probability measures with a fixed finite support.
A subset of $\G$, say $S$, is ``symmetric'' if $x^{-1}\in S$ whenever
$x\in S$.
Let $S$ be a finite symmetric subset of $\G$. We assume that $S$
generates the whole group $\G$.
Let ${\cal P}_s(S)$ be the set of symmetric probability measures with support
equal to $S$. Then ${\cal P}_s(S)$ naturally identifies with an open
subset of $\R^d$ for some $d$. We use the differentiable structure it
inherits this way.

\emph{Regularity Assumption}: throughout the paper, we shall assume that
the function $\l\in[-1,1]\ra\log\mu_\l(a)$ has a derivative
at $\l=0$ for all $a\in S$. Equivalently, we may write a first order
expansion of $\log\mu_\l(a)$ as $\l$ tends to $0$ in the form
%
%e2.3 #&#
\begin{equation}
\label{eq:reg}\log\mu_\l(a)=\log\mu_0(a)+\l\nu (a)+\l
o_\l(a),
\end{equation}
where $\nu(a)$ is the derivative of the function $\l\ra\log\mu_\l(a)$
at $\l=0$ and
$o_\l(a)$ converges to $0$.

Observe that since $S$ is finite, this is equivalent to requiring $o_\l
(a)$ to converge to $0$ uniformly with respect to $a\in S$.
We shall also repeatedly use the fact that $\nu$ is bounded.

We shall use the shorthand notation $\P^\l$ (resp., $\esp^\l$) instead
of $\P^{\mu_\l}$ (resp., $\esp^{\mu_\l}$).

From the condition that $\mu_\l$ is a probability measure, one deduces
that we must have
\[
\sum_{a\in S} \nu(a)\mu_0(a)=0.
\]
We define the sequence $M_0=0$ and, for $n\geq1$,
\[
M_n=\sum_{j=1}^n
\nu(X_j).
\]
Note that the random process $(M_n)_{n\in\N}$ is a centered
martingale under $\P^0$.

Let $d\in\D(\G)$ and assume assumption (\emph{BA}) holds.
We shall see in Proposition~\ref{prop:clt} that the sequence $(\vert
Z_n\vert,M_n)$ satisfies a joint central limit theorem and that the asymptotic
covariance of $\vert Z_n\vert$ and $M_n$ is given by
%
%e2.4 #&#
\begin{equation}
\label{eq:cov}\sigma(\nu,\mu_0;d):=\lim_{n\ra
\infty}
\frac{1} n \esp ^0\bigl[\vert Z_n\vert M_n\bigr].
\end{equation}

%th2.1 #&#
\begin{theo}
\label{theo:rescp}
Let $d\in\D(\G)$ satisfy {(BA)}. Then the function $\l\ra
\ell(\mu_\l
;d)$ is differentiable and its derivative satisfies
%
%e2.5 #&#
\begin{equation}
\label{eq:rescp} \frac{d}{d\l}\bigg\vert _{\l=0} \ell(
\mu_\l;d) = \sigma(\nu,\mu_0;d).
\end{equation}
\end{theo}

%th2.2 #&#
\begin{theo}
\label{theo:entrop}
The function $\l\ra h(\mu_\l)$ has a derivative at $\l=0$ which satisfies
%
%e2.6 #&#
\begin{equation}
\label{eq:entrop} \frac{d}{d\l}\bigg\vert _{\l=0} h(
\mu_\l) = \sigma\bigl(\nu,\mu_0;d_G^0
\bigr),
\end{equation}
where $d_G^\l:=d_G^{\mu_\l}$ is the Green metric associated to the
probability $\mu_\l$.
\end{theo}

In \cite{kn:bhm1}, we showed that the asymptotic entropy coincides with
the rate of escape in the Green metric:
\[
h(\mu)=\ell\bigl(\mu;d_G^\mu\bigr)
\]
for all $\mu\in{\cal P}_s(S)$.
Thus, using the result in Theorem~\ref{theo:rescp}, we may reformulate
(\ref{eq:entrop}) as follows:
\[
\lim_{\l\ra0}\frac{1}\l\bigl( h(\mu_\l)-h(
\mu_0)\bigr)=\lim_{\l\ra
0}\frac{1}\l\bigl( \ell
\bigl(\mu_\l;d_G^0\bigr)-\ell\bigl(
\mu_0;d_G^0\bigr)\bigr).
\]
In other words, as far as first-order terms are concerned, the
fluctuations of the Green metric do not contribute to the fluctuations
of the entropy. This turns out to be a quite general statement for
random walks on nonamenable groups; see Section~\ref{sec:entrop}.

%%%%%%%%%%%%%%%%%%%%%%%%%%%%%%%%%%%%%%%%%%%%%%%%%%%%%%%%%%%%%%
%s2.3 #&#
\subsection{Heuristics} \label{sec:heu}

We give a simple---but not completely rigorous---way to guess why
formula (\ref{eq:rescp})
%and (\ref{eq:entrop})
should hold true. We provide these heuristics in order to clarify the
scheme of the proofs, with the hope that this scheme can be adapted to
other examples of random walks.

To compute the rate of escape, observe that
%
%e2.7 #&#
\begin{equation}
\label{eq:gir} \esp^\l\bigl[\vert Z_n\vert\bigr]=\esp^0
\Biggl[\vert Z_n\vert\prod_{j=1}^n
\frac
{\mu_\l
(X_j)}{\mu_0(X_j)}\Biggr].
\end{equation}

Taking the derivative in (\ref{eq:gir}), we get that
%
%e2.8 #&#
\begin{equation}
\label{eq:gir2} \frac{d}{d\l}\bigg\vert _{\l=0}
\esp^\l\bigl[\vert Z_n\vert\bigr]=\esp ^0\Biggl[\vert
Z_n\vert\sum_{j=1}^n
\nu(X_j)\Biggr] =\esp^0\bigl[\vert Z_n\vert
M_n\bigr].
\end{equation}

Thus, we see that a reasonable candidate to be the derivative of $\ell
(\mu_\l;d)$ is the limit of
$\frac{1} n \esp^0[\vert Z_n\vert M_n]$ as $n$ tends to $+\infty$.

Observe, however, that in order to turn this loose argument into a
proof, one needs justify how to exchange the order between the limit in
$n$ and the derivation in $\l$.

On the one hand, we shall rely on a quantitative version of the law of
large numbers for $\vert Z_n\vert$ to show that the derivative of
$\ell
(\mu_\l;d)$ is well approximated by the limit of the ratio $(\esp^\l
[\vert Z_n\vert]-\esp^0[\vert Z_n\vert])/{\l n}$ as soon as $\l$ tend
to $0+$ and $n$ tend to $+\infty$ in such a way that the product $\l n$
tends to $+\infty$, see Lemma~\ref{lem:geoapp}. This follows from the
fact that the function $n\ra\esp^\l[\vert Z_n\vert]$ is almost
additive, uniformly in $\l$.

Thus, it is sufficient to describe the limit of $\esp^\l[\vert
Z_n\vert
]/\l n$ when $\l n\ra+\infty$. We will actually choose $\l$ and $n$
such that $\l\sqrt{n}$ tends to $1$. Then $\esp^\l[\vert Z_n\vert
]/\l
n\sim\esp^\l[\vert Z_n\vert]/\sqrt n$, which is the scaling of the
central limit theorem.

More precisely, we show a joint C.L.T. for the random vector $(\vert
Z_n\vert,M_n)$ under~$\P^0$; see Proposition~\ref{prop:clt}. Let
$\sigma
(\nu,\mu_0;d)$ denote the asymptotic covariance of $\vert Z_n\vert$
and $M_n$.

Consider the Girsanov weight $\prod_{j=1}^n {\mu_\l(X_j)}/{\mu_0(X_j)}$
in formula (\ref{eq:gir}).
Up to error terms of smaller order, it coincides with the exponential
(in the sense of martingale theory) of the martingale $(\l M_n)_{n\in
\N
}$. With our choice of the scaling between $\l$ and $n$, the
asymptotic of
$\l M_n\sim M_n/\sqrt n$ is given by the central limit theorem.
Therefore, the limit of the Girsanov weight is of the form
 $e^{M-({1} /{2})\esp[M^2]}$ for some Gaussian random variable $M$.
Moreover, provided we can check some integrability conditions, the
joint C.L.T. implies that the limit of
$\esp^\l[\vert Z_n\vert]-n\ell(\mu_0;d)]/\l n$
is then of the form $\esp[Ze^{M-({1}/ {2})\esp[M^2]}]$, where $(Z,M)$
is a Gaussian vector with covariance
$\esp[ZM]=\sigma(\nu,\mu_0;d)$. The integration by parts formula for
Gaussian laws implies that (for any Gaussian vector)
\[
\esp\bigl[Ze^{M-({1}/ 2)\esp[M^2]}\bigr] =\esp[ZM].
\]

The next theorem summarizes the part of this argument we just sketched
that does not explicitly use the hyperbolicity of $\G$.

%th2.3 #&#
\begin{theo} \label{theo:malliavin}
Let $\G$ be a finitely generated group; $S$ a finite symmetric
generating set; $d$ a left-invariant proper metric on $\G$ and $\l\in
[-1,1]\ra\mu_\l\in{\cal P}_s(S)$ be a curve in ${\cal P}_s(S)$
satisfying the
Regularity Assumption. We further assume that:
\begin{longlist}[(ii)]
\item[(i)] the joint central limit theorem holds for the vector $(\vert
Z_n\vert, M_n)$ under $\P^0$ with asymptotic covariance $\sigma$,\vspace*{3pt}
%%and that $\frac1 n \esp^0[\vert Z_n\vert M_n]$ converges to $\sigma$
%as $n$ tends to $+\infty$,\\{\chr"CA}

\item[(ii)]
\[
\sup_n \frac{1} n\esp^0\bigl[\bigl(\vert
Z_n\vert-n\ell(\mu_0;d)\bigr)^2\bigr]<+
\infty.
\]
Then we have
\[
\lim_{n\rightarrow+\infty,\l\rightarrow0} \frac{1}{\l n} \bigl(\esp^\l \bigl[\vert
Z_n\vert\bigr]-\esp^0\bigl[\vert Z_n\vert\bigr]\bigr)=
\sigma,
\]
along any sequence $\l$ such that $\limsup_{n\rightarrow+\infty} \l^2
n<+\infty$.
\end{longlist}
\end{theo}

Theorem~\ref{theo:malliavin} is proved in Section~\ref{sec:malliavin}.

%Let us now turn to the entropy.
%
%We start with the definition $$h(\mu_\l)=\lim_n \frac1 n \sum_{x\in
%and observe that
%$$\sum_{x\in\G} \left.\frac d{d\l}\right|_{\l=0} (-\log\mu_\l^n(x))
%=\sum_{x\in\G} -\left.\frac d{d\l}\right|_{\l=0} \mu_\l^n(x)
%=0$$
%since $\sum_{x\in\G} \mu_\l^n(x)=1$ for all $\l$.
%
%Assuming that we may exchange limits, we would then get that
%$$\left.\frac d{d\l}\right|_{\l=0} h(\mu_\l)
%= \left.\frac d{d\l}\right|_{\l=0} \lim_n \frac1 n \sum_{x\in\G} (-
%=\left.\frac d{d\l}\right|_{\l=0} \ell(\mu_\l;d_G^0)
%.$$
%In other words we would have proved that the derivative of the entropy
%and the derivative of the rate of escape in the Green metric coincide.
%We refer to part \ref{sec:entrop} for a rigorous statement.
%
%%%%%%%%%%%%%%%%%%%%%%%%%%%%%%%%%%%%%%%%%%%%%%%%%%%%%%%%%%%%%%%%%%%%%%%%%%%%%%%%%%%%%
%%%%%%%%%%%%%%%%%%%%%%%%%%%%% PROOF THEOREM RESCP
%%%%%%%%%%%%%%%%%%%%%%%%%%%%%%%%%%%%%%%%%%%%%%%
%%%%%%%%%%%%%%%%%%%%%%%%%%%%%%%%%%%%%%%%%%%%%%%%%%%%%%%%%%%%%%%%%%%%%%%%%%%%%%%%%%%%%

%s3 #&#
\section{Proofs of Theorems \texorpdfstring{\protect\ref{theo:rescp}}{2.1} and \texorpdfstring{\protect\ref{theo:malliavin}}{2.3}}
\label{sec:proof}

The proofs are organized in the following way: in Section \ref
{sec:geometricinput}, we recall some estimates on the mean distance
$\esp^\mu[\vert Z_n\vert]$ from~\cite{kn:bhm2} and use them to show
that Theorem~\ref{theo:rescp} can be deduced from Theorem~\ref
{theo:malliavin}. In Section~\ref{sec:centrallimittheorems}, we recall
results from \cite{kn:bjork} and show how they imply a slightly
stronger version of the assumptions of Theorem~\ref{theo:malliavin}.
Sections~\ref{sec:geometricinput} and \ref{sec:centrallimittheorems} use
the hyperbolicity of $\G$ in an essential way.

In Section \ref{sec:malliavin}, we prove Theorem~\ref{theo:malliavin}.
Section \ref{sec:malliavin} can be read independently of the preceding ones.

Let $d\in\D(\G)$ satisfy assumption \emph{(BA)}.
We will show that\break $\lim_{\l\ra0}\frac{1}\l(\ell(\mu_\l;d)-\ell
(\mu
_0;d))$ exists and equals $\sigma(\nu,\mu_0;d)$. In the proof, it will
be convenient to restrict ourselves to positive $\l$'s. This is no loss
of generality since $\sigma(\nu,\mu_0;d)$ is linear in $\nu$.

%s3.1 #&#
\subsection{Geometric input}
\label{sec:geometricinput}

%le3.1 #&#
\begin{lm} \label{lem:geoapp} Let $\l$ tend to $0+$ and $n$ tend to
$+\infty$ in such a way that the product $\l n$ tends to $+\infty$. Then
\[
\frac{\ell(\mu_\l;d)-\ell(\mu_0;d)}\l-\frac{\esp^\l[\vert
Z_n\vert]-\esp
^0[\vert Z_n\vert]}{\l n}
\]
tends to $0$.
\end{lm}

\begin{pf}
Let $\mu\in{\cal P}_s(S)$.

The triangle inequality implies that the sequence $a(n):=\esp^\mu
[\vert
Z_n\vert]$ is sub-additive. Therefore, we have $\ell(\mu;d)=\inf_n\frac
{a(n)}n$, and thus
%
%e3.1 #&#
\begin{equation}
\label{eq:dev0} \esp^\mu\bigl[\vert Z_n\vert\bigr]\geq n\ell (\mu;d)
\end{equation}
for all $n$.

We need a similar upper bound. It will follow from bounds on the
lateral deviation of a trajectory of a random walk. Let us recall some
results from \cite{kn:bhm2}.

In \cite{kn:bhm2}, Proposition~3.8, we showed that, for any $\mu\in
{\cal P}
_s(S)$, there exists a constant $\tau_0$ such that for all integers
$m$, $n$, $k$,
%
%e3.2 #&#
\begin{equation}
\label{eq:dev2}\esp^\mu \bigl[(Z_m,Z_{m+n+k})_{Z_{m+n}}
\bigr]\leq\tau_0.
\end{equation}

Applying this inequality with $m=0$ and using the fact that $\esp^\mu
[\vert Z_{n+k}-Z_n\vert]=\esp^\mu[\vert Z_k\vert]$, we get that
%
%e3.3 #&#
\begin{equation}
\label{eq:eq:dev3} \esp^\mu\bigl[\vert Z_n\vert\bigr]+\esp
^\mu\bigl[\vert Z_k\vert \bigr]\leq2\tau_0+
\esp^\mu\bigl[\vert Z_{n+k}\vert\bigr].
\end{equation}
Thus, the sequence $b(n):=2\tau_0-\esp^\mu[\vert Z_n\vert]$ is also
subadditive. Note that $b(n)/n$ converges to $-\ell(\mu;d)$. As above,
this implies that
%
%e3.4 #&#
\begin{equation}
\label{eq:dev4} \esp^\mu\bigl[\vert Z_n\vert\bigr]\leq n\ell (\mu
;d)+2\tau _0
\end{equation}
for all $n$.

Combining (\ref{eq:dev0}) and (\ref{eq:dev4}), we see that we have
proved that
%
%e3.5 #&#
\begin{equation}
\label{eq:dev} \bigl\vert\esp^\mu\bigl[\vert Z_n\vert\bigr]- n\ell(\mu
;d)\bigr\vert \leq2\tau_0
\end{equation}
for all $n$.

A close inspection of the proof of (\ref{eq:dev2}) reveals that the
constant $\tau_0$ is locally uniform in ${\cal P}_s(S)$ so that we may apply
(\ref{eq:dev}) with the same constant $\tau_0$ to all measures $\mu
_\l$
for $\l$ in a small enough neighborhood of $0$.
The statement of Lemma~\ref{lem:geoapp} immediately follows.
\end{pf}

Let us choose $\l$ tending to $0+$ and $n$ tending to $+\infty$ such
that $\l^2n$ tends to $1$.
Thus, Lemma~\ref{lem:geoapp} applies. In order to complete the proof of
Theorem~\ref{theo:rescp},
it only remains to show that
%
%e3.6 #&#
\begin{equation}
\label{eq:fond} \lim_{\l,n} \frac{\esp^\l[\vert Z_n\vert]-\esp^0[\vert Z_n\vert
]}{\l
n}=\sigma(\nu,
\mu_0;d).
\end{equation}

%%%%%%%%%%%%%%%%%%%%%%%%%%%%%%%%%%%%%%%%%%%%%%%%%%%%%%%%%%%%%%%%%%%%%%%%%%%%
%s3.2 #&#
\subsection{Central limit theorems}
\label{sec:centrallimittheorems}

In this part of the paper, we recall some results from \cite{kn:bjork}
on the central limit theorem for $\vert Z_n\vert$ and discuss their
extension to a joint C.L.T. for $(\vert Z_n\vert,M_n)$.

Let $d\in\D(\G)$ satisfy assumption (\emph{BA}).
Let $\mu\in{\cal P}_s(S)$ be a finitely supported symmetric probability
measure on $\G$.

Mimicking the situation of the discussion preceding Theorem~\ref
{theo:rescp}, we also let $\nu$ be a real valued function defined on
$S$ and satisfying the centering condition: $\sum_{a\in S} \nu(a)\mu
(a)=0$ and
consider the sequence of random variables $M_0=0$ and, for $n\geq1$,
\[
M_n=\sum_{j=1}^n
\nu(X_j).
\]

%pr3.2 #&#
\begin{prop}\label{prop:clt}
\textup{(i)} The law of the two-dimensional random vector
$( (\vert Z_n\vert-n\ell(\mu;d))/\sqrt{n},M_n/\sqrt{n})$ under $\P
^\mu$
weakly converges as $n$ tends to $+\infty$ to a centered Gaussian law
with some covariance matrix $\Sigma^\mu$.

\textup{(ii)} The covariance matrix of $( (\vert Z_n\vert-n\ell(\mu;d))/\sqrt {n},M_n/\sqrt{n})$ under $\P^\mu$ converges to $\Sigma^\mu$. \\
In particular, the sequence $\frac{1} n \esp^\mu[\vert Z_n\vert M_n]$
converges as $n$ tends to $+\infty$ and its limit is the nondiagonal
term of $\Sigma^\mu$.
\end{prop}

\begin{pf}
We recall the following classical version of the martingale central
limit theorem (see \cite{kn:hel}).

%le3.3 #&#
\begin{lm} \label{lem:clt}
Let $(\zeta_n)_{n\in\N}$ be a square integrable, centered martingale
with respect to a filtration
$({\cal F}_n)_{n\in\N}$, with stationary increments.
Assume that
%
%e3.7 #&#
\begin{equation}
\label{eq:lind}\frac{1} n \sum_{j=1}^n
\esp\bigl[(\zeta _j-\zeta _{j-1})^2\vert{\cal
F}_{j-1}\bigr]\ra\sigma^2
\end{equation}
almost surely, where $\sigma^2$ is a deterministic real. Then the law
of $(\frac{1} {\sqrt{n}} \zeta_n)$ weakly converges to a centered
Gaussian law with variance $\sigma^2$.
\end{lm}

\textit{Step} 1:
Following \cite{kn:bjork}, we first prove a version of Proposition~\ref
{prop:clt} where $\vert Z_n\vert-n\ell(\mu;d)$ is replaced by an
appropriate martingale approximation that we denote with $(\chi_n)$.

Let $\partial\G$ the Busemann boundary of $\G$. We recall that the
Gromov product extends to the boundary.

Under $\P^\mu$, almost any trajectory $(Z_n)_{n\in N}$ converges to a
limiting point in $\partial\Gamma$, say $Z_\infty$.
This follows from assumption (\emph{BA}) since then $\partial\G$
can be
identified with the Gromov boundary of $\G$ and one knows that random
walk paths almost surely converge in the Gromov compactification of a
hyperbolic group; see \cite{kn:anc} or \cite{kn:kaim}.
The law of $Z_\infty$ is called the ``harmonic measure.'' We denote it
with $\xi^\mu$.

In \cite{kn:bjork}, part 4, it is proved that there exists a bounded
function $u$ on $\partial\G$ such that
the sequence
\[
\chi_n:=k(Z_n)-n\ell(\mu;d)+u(k)-u\bigl(Z_n^{-1}k
\bigr)
\]
is a martingale under $\P^\mu$ for any $k\in\partial\G$.

In the sequel, we shall assume that $k$ is chosen according to the
harmonic measure $\xi^\mu$ and independent of the walk $(Z_n)_{n\in
\N}$.
It then follows that $(\chi_n)_{n\in\N}$ has stationary increments.
It is also proved in Theorem~9 in \cite{kn:bjork} that the Lindeberg
condition (\ref{eq:lind}) is satisfied.

%re3.4 #&#
\begin{rmk}
The group $\G$ has a natural action on its boundary. For each $k\in
\partial\G$, the sequence $(Z_n^{-1}k)_{n\in\N}$ is a Markov chain
with values in $\partial\G$ started at $k$. This Markov chain has a
unique invariant probability measure, namely the harmonic measure $\xi
^\mu$. Moreover, the Markov chain $(Z_n^{-1}k)_{n\in\N}$ has nice
mixing properties; see Lemma~4 in \cite{kn:bjork}. By choosing $k$
according to $\xi^\mu$, we are actually considering the chain in a
stationary regime.

We shall not explicitly use these remarks but they lie at the heart of
the proof of the central limit theorem in \cite{kn:bjork}.
\end{rmk}

On the other hand, the sequence $(M_n)_{n\in\N}$ being a sum of
independent bounded and centered random variables, is also a centered
martingale with stationary increments
under $\P^\mu$ satisfying condition (\ref{eq:lind}).

Thus, we deduce from Lemma~\ref{lem:clt} that the law of the vector
$(\chi_n/\sqrt{n},\break M_n/\sqrt{n})$ under $\P^\mu$ converges to a centered
Gaussian vector. Indeed, one may apply Lemma~\ref{lem:clt} to the
martingale $(a\chi_n+bM_n)_{n\in\N}$ for any $a,b\in\R$.

Let $\Sigma^\mu$ be the limit covariance matrix. For $(a,b)\in\R
^2$, we
use the notation $\Sigma^\mu(a,b)$ to denote the value of the quadratic
form associated to $\Sigma^\mu$ evaluated at $(a,b)$.
We observe that since both martingales $(\chi_n)$ and $(M_n)$ have
stationary increments, then $\Sigma^\mu$ is also the covariance (under
$\P^\mu$) of the vector $(\chi_n/\sqrt{n},M_n/\sqrt{n})$ for any
$n\geq
1$, that is,
%
%e3.8 #&#
\begin{equation}
\label{eq:bd1} \frac{1} n\esp^\mu\bigl[(a\chi_n+bM_n)^2
\bigr]=\Sigma^\mu(a,b).
\end{equation}

\textit{Step} 2:
In the above claims, we wish to replace $\chi_n$ by $\vert Z_n\vert
-n\ell(\mu;d)$.
We shall use the following lemma.

%le3.5 #&#
\begin{lm}\label{lem:shad}
There exists a constant $C$ such that, for all $D$ we have
%
%e3.9 #&#
\begin{equation}
\label{eq:bd4} \P^\mu\bigl[(k,Z_n)_{\mathrm{id}}\geq D
\bigr]\leq C^{-1}e^{-CD}.
\end{equation}
Inequality (\ref{eq:bd4}) holds uniformly in $k\in\partial\G$.
\end{lm}

\begin{pf}%{Proof of Lemma~\ref{lem:shad}}
The statement of the lemma actually directly follows from arguments in
\cite{kn:bhm2}.

One may for instance split the event $(k,Z_n)_{\mathrm{id}}\geq D$ into two, say
$A:=((k,Z_n)_{\mathrm{id}}\geq D)\cap((Z_\infty,Z_n)_{\mathrm{id}}\geq\frac{D} 2)$ and
$B:=((k,Z_n)_{\mathrm{id}}\geq D)\cap((Z_\infty,Z_n)_{\mathrm{id}}< \frac{D} 2)$.

Let us show that $\P^\mu[A]+\P^\mu[B]\leq C^{-1}e^{-CD}$.

In the argument below, $\tau_1$ is a constant that depends on $d$ and
$\mu$ only. We choose $D$ large enough and how large depends only on
the choice of the metric $d$ and the measure $\mu$. In particular, we
assume that $D\geq4\tau$, where $\tau$ is the hyperbolicity constant
from (\ref{eq:hyper}).

Hyperbolicity implies that, on $A$, we also have $(k,Z_\infty
)_{\mathrm{id}}\geq
\frac{D} 2-\tau_1\geq\frac{D} 2$.
We know from \cite{kn:bhm2}, Proposition~3.10, that $\xi^\mu$ satisfies
the doubling condition.
Therefore, the probability that $(k,Z_\infty)_{\mathrm{id}}\geq\frac{D} 2$ can
be compared to the harmonic measure of a ball
of $\partial\Gamma$ of radius of order $e^{-C_1 D}$ for some
$C_1$, and since $\xi^\mu$ is Ahlfors regular (see Theorem~1.1 in
\cite
{kn:bhm2}) we get that
$\P^\mu[A]\leq C^{-1}e^{-CD}$.

On the event $B$, we have $(Z_\infty,Z_n)_{\mathrm{id}}< \frac{D} 2$ and $\vert
Z_n\vert\geq D-\tau\geq\frac{3}4 D$. Therefore, the distance between
$Z_n$ and any
quasiruler from $\mathrm{id}$ to $Z_\infty$ is larger than $\frac{D} 4-\tau_1$.
For large enough $D$, this last event has a probability bounded from above
by $C^{-1}e^{-CD}$ for some $C$ as follows from the deviation
inequality in Proposition~3.8 in \cite{kn:bhm2}. Therefore,
$\P^\mu[B]\leq C^{-1}e^{-CD}$.
\end{pf}

Back to the proof of Proposition~\ref{prop:clt},
we observe that $k(x)=\vert x\vert-2(k,x)_{\mathrm{id}}$ for all $k\in\partial
\G
$ and $x\in\G$.
Therefore,
\[
\vert Z_n\vert-n\ell(\mu;d)-\chi _n=2(k,Z_n)_{\mathrm{id}}-
\bigl(u(k)-u\bigl(Z_n^{-1}k\bigr)\bigr).
\]
Using Lemma~\ref{lem:shad} and the fact that $u$ is bounded, we get that
%
%e3.10 #&#
\begin{equation}
\label{eq:bd3}\P^\mu\bigl[\bigl\vert\vert Z_n\vert-n\ell (
\mu;d)-\chi_n \bigr\vert\geq D\bigr]\leq C^{-1}e^{-CD}
\end{equation}
for some constant $C$.

As a by-product of (\ref{eq:bd3}), we get that $\frac{1}{\sqrt{n}}
(\vert Z_n\vert-n\ell(\mu;d)-\chi_n)$ converges to $0$
in probability. Therefore, the two sequences of vectors $( (\vert
Z_n\vert-n\ell(\mu;d))/\sqrt{n},M_n/\sqrt{n})$ and
$( \chi_n/\sqrt{n},M_n/\sqrt{n})$ have the same limit in law. Thus, we
have proved that $( (\vert Z_n\vert-n\ell(\mu;d))/\sqrt {n},M_n/\sqrt {n})$ satisfies the central limit theorem with asymptotic variance
$\Sigma^\mu$.

We also deduce from (\ref{eq:bd3}) that
%
%e3.11 #&#
\begin{equation}
\label{eq:bd2}\sup_n \esp^\mu\bigl[\bigl(\vert
Z_n\vert -n\ell(\mu;d)-\chi _n\bigr)^2
\bigr]<\infty
\end{equation}
and
\[
\frac{1} n \esp^\mu\bigl[\bigl(\vert Z_n\vert-n
\ell(\mu;d)-\chi_n\bigr)^2\bigr]\ra0.
\]

Combining (\ref{eq:bd2}) with (\ref{eq:bd1}), we see that the
covariance matrix of $( (\vert Z_n\vert-n\ell(\mu;d))/\sqrt {n},M_n/\sqrt {n})$ under $\P^\mu$ converges to $\Sigma^\mu$.
Also
\[
\frac{1} n \esp^\mu\bigl[\vert Z_n\vert
M_n\bigr]=\frac{1} n \esp^\mu\bigl[\bigl(\vert
Z_n\vert-n\ell(\mu;d)\bigr) M_n\bigr]
\]
converges as $n$ tends to $+\infty$ and its limit is the nondiagonal
term of $\Sigma^\mu$.
This concludes the proof of the proposition.
\end{pf}

For further references, we observe that (\ref{eq:bd2}) with (\ref
{eq:bd1}) implies the following.

%le3.6 #&#
\begin{lm}\label{lem:2ndmoment}
For all $\mu\in{\cal P}_s(S)$, we have
\[
\sup_n \frac{1} n\esp^\mu\bigl[\bigl(\vert
Z_n\vert-n\ell (\mu;d)\bigr)^2\bigr]<\infty.
\]
\end{lm}

%%%%%%%%%%%%%%%%%%%%%%%%%%%%%%%%%%%%%%%%%%%%%%%%%%%%%%%%%%%%%%%%%%%%%%%
%s3.3 #&#
\subsection{Proof of Theorem \texorpdfstring{\protect\ref{theo:malliavin}}{2.3} and (\texorpdfstring{\protect\ref{eq:fond}}{3.6})}
\label{sec:malliavin}

Here, $\G$ is any finitely generated group.

We recall that we are assuming the joint central limit theorem for the
vector $(\vert Z_n\vert, M_n)$ under $\P^0$; namely the law of the
two-dimensional random vector
$( (\vert Z_n\vert-n\ell(\mu;d))/\sqrt{n},M_n/\sqrt{n})$ under $\P^0$
weakly converges as $n$ tends to $+\infty$ to a centered Gaussian law
with some covariance matrix $\Sigma$. Let $\sigma$ be the nondiagonal
element of $\Sigma$.

We also assume that
%
%e3.12 #&#
\begin{equation}
\label{eq:varass}\sup_n \frac{1} n\esp^0
\bigl[\bigl(\vert Z_n\vert-n\ell (\mu _0;d)
\bigr)^2\bigr]<+\infty.
\end{equation}

We wish to prove that
%
%e3.13 #&#
\begin{equation}
\label{eq:fond'}\lim_{n\rightarrow+\infty,\l
\rightarrow0} \frac{1}{\l n} \bigl(
\esp^\l\bigl[\vert Z_n\vert\bigr]-\esp^0\bigl[\vert
Z_n\vert\bigr]\bigr)=\sigma,
\end{equation}
along any sequence $\l$ such that $\limsup_{n\rightarrow+\infty} \l^2
n<+\infty$.
Without loss of generality, we may and will assume that $\l^2n$
converges to some limit $\alpha\ge0$. Then
$\l\sim\sqrt{\alpha/n}$.

We start dealing with the case $\alpha\neq0$.

First, note that
\[
\frac{\esp^\l[\vert Z_n\vert]-\esp^0[\vert
Z_n\vert]}{\l n} = \frac{\esp^\l[\vert Z_n\vert-n\ell(\mu_0;d)]-\esp^0[\vert
Z_n\vert
-n\ell(\mu_0;d)]}{\l n}.
\]

From the central limit theorem for $\vert Z_n\vert$ under $\esp^0$ and
assumption (\ref{eq:varass}), we get that
\[
\frac{\esp^0[\vert Z_n\vert-n\ell(\mu_0;d)]}{\l n} \sim\frac{\esp^0[\vert Z_n\vert-n\ell(\mu_0;d)]}{\sqrt{\alpha
}\sqrt {n}}\ra0.
\]

In order to compute the limit of $\esp^\l[\vert Z_n\vert-n\ell(\mu
_0;d)]/\l n$, we write these terms in a form
that is more amenable to the application of the central limit theorem.

We have
%
%e3.14 #&#
\begin{equation}\qquad\qquad
\label{eq:fond0} \frac{1}{\l n} \esp^\l\bigl[\vert Z_n
\vert-n\ell(\mu_0;d)\bigr]=\frac{1}{\l n} \esp ^0
\Biggl[\bigl(\vert Z_n\vert-n\ell(\mu_0;d)\bigr)
\prod_{j=1}^n \frac{\mu_\l
(X_j)}{\mu
_0(X_j)}\Biggr].
\end{equation}

Let $a\in S$. Recall that we have a first-order expansion of the
function $\l\ra\mu_\l(a)$ in the form
%
%e3.15 #&#
\begin{equation}
\label{eq:exp1} \log\frac{\mu_\l(a)}{\mu_0(a)}= \l\nu(a)+\l o_\l(a),
\end{equation}
where the function $o_\l$ uniformly converges to $0$ as $\l$ goes to $0$.

Because $\mu_\l$ is a probability for all $\l$, it follows from
(\ref
{eq:exp1}) that $\nu$ and $o_\l$ must
satisfy the following centering conditions:
%
%e3.16 #&#
%e3.17 #&#
\begin{eqnarray}
\label{eq:center0} \sum_{a\in S} \nu(a)
\mu_0(a)&=&0\quad\mbox{and }
\\
\label{eq:center}  \lim_{\l\ra0}\frac{1}\l
\sum_{a\in S} \biggl(o_\l(a)+\frac\l2
\nu^2(a)\biggr)\mu_0(a)&=&0.
\end{eqnarray}

To see why (\ref{eq:center0}) and (\ref{eq:center}) hold, note that
$\sum_{a\in S} \mu_\l(a)=1$ for all $\l$.
The expansion of $\sum_{a\in S} \mu_\l(a)$ in terms of $\l$ starts with
$\l\sum_{a\in S} \nu(a)\mu_0(a) +\l^2\frac{1}\l\sum_{a\in S}
(o_\l
(a)+\frac\l2 \nu^2(a))\mu_0(a)$, the rest being of order smaller than
$\l^2$.
Dividing by $\l$ and letting $\l$ tend to $0$, one gets (\ref
{eq:center0}). Then dividing by $\l^2$ and letting $\l$ tend to $0$,
one gets (\ref{eq:center}).

Let us rewrite (\ref{eq:fond0}) as
%
%e3.18 #&#
\begin{equation}\qquad\qquad
\label{eq:fond1} \frac{1}{\l n} \esp^\l\bigl[\vert Z_n
\vert-n\ell(\mu_0;d)\bigr]=\frac{1}{\l n} \esp ^0\bigl[
\bigl(\vert Z_n\vert-n\ell(\mu_0;d)\bigr) e^{\l M_n-\l^2 A_n+ R^\l_n}
\bigr],
\end{equation}
where
\[
M_n:=\sum_{j=1}^n
\nu(X_j), \qquad\hspace*{-2pt} A_n:=\frac{1}2\sum
_{j=1}^n \nu^2(X_j) \hspace*{-1pt}\quad\mbox{and}\hspace*{-1pt} \quad R^\l_n=\l\sum_{j=1}^n
o_\l(X_j)+\frac\l2\nu^2(X_j).
\]

From the Law of Large Numbers for sums of i.i.d. random variables, it
follows that $\frac{1} n A_n$ almost surely converges
to $\frac{1}2 \sum_{a\in S} \nu^2(a)\mu_0(a)$ and, therefore, $\l^2A_n$
converges to $\frac{\alpha}2 \sum_{a\in S} \nu^2(a)\mu_0(a)$. We claim
that $R^\l_n$ converges to $0$ in probability under $\P^0$.

The argument for this last claim goes as follows. Let $Y^\l_j:=\frac{1}\l o_\l(X_j)+\frac{1} 2\nu^2(X_j)$, so that
$R^\l_n=\l^2\sum_{j=1}^nY^\l_j$.
For a fixed $n$, the random variables $(Y^\l_j)_{j=1}^n$ are
independent and equally distributed.
On the one hand, (\ref{eq:center}) implies that $\esp^0[Y^\l_1]$ tends
to $0$ and, therefore, $\esp^0[R^\l_n]$ also converges to $0$. On the
other hand, the variance of $Y^\l_j$ satisfies $\lim_{\l\ra0} \l
^2\var
^0[Y^\l_1]=0$. (Use the fact that $o_\l$ converges to $0$.)
Therefore, the variance of $R^\l_n$ is of lower order than
$\l^2n$ and tends to $0$. Thus, we get that both the mean and variance
of $R^\l_n$ converge to $0$.

We now use the central limit theorem for
$((\vert Z_n\vert-n\ell(\mu_0;d))/\sqrt{n}, M_n/\sqrt{n})$.

Ignoring for a moment that the function we integrate in (\ref
{eq:fond1}) is not bounded, we pass to the limit using the relation $\l
^2 n\ra\alpha$ and get that
%
%e3.19 #&#
\begin{equation}
\label{eq:fond2} \frac{1}{\l n} \esp^\l\bigl[\vert Z_n
\vert-n\ell(\mu_0;d)\bigr] \ra\frac{1}{\sqrt\alpha} \esp
\bigl[Ze^{\sqrt{\alpha}M-(\alpha/ 2)\sum
_{a\in S} \nu^2(a)\mu_0(a)}\bigr],
\end{equation}
where $(Z,M)$ is a centered Gaussian vector with variance $\Sigma$.

Observe that the variance of $M_n$ equals $n\sum_{a\in S} \nu^2(a)\mu
_0(a)$ and, therefore,
$\sum_{a\in S} \nu^2(a)\mu_0(a)=\esp[M^2]$. Thus, the right-hand side
of (\ref{eq:fond2}) equals
\[
\frac{1}{\sqrt\alpha}\esp\bigl[Ze^{\sqrt{\alpha}M-(\alpha/2)\esp
[M^2]}\bigr].
\]
The integration by parts formula for Gaussian laws implies that (for
any Gaussian vector and any $\alpha$)
\[
\frac{1}{\sqrt\alpha}\esp\bigl[Ze^{\sqrt{\alpha}M-(\alpha/2)\esp[M^2]}\bigr] =\esp[ZM]=\sigma.
\]

Thus, we are done with the proof of Theorem~\ref{theo:rescp} once we
justify that we may indeed pass to the limit in (\ref{eq:fond1}).
In order to do so, it is sufficient to have bounds on the moments of
the functions we integrate.

H\"older's inequality implies that
\begin{eqnarray*}
&&\esp^0\biggl[\biggl(\frac{1}{\l n}\bigl(\vert Z_n
\vert-n\ell(\mu _0;d)\bigr) e^{\l
M_n-\l^2 A_n+ R^\l_n}\biggr)^{6/5}
\biggr]
\\
&&\qquad\leq \esp^0\biggl[\frac{1}{(\l n)^2}\bigl(\vert Z_n
\vert-n\ell(\mu_0;d)\bigr)^2\biggr]^{6/10}
\esp^0\bigl[e^{3\l M_n-3\l^2 A_n+3 R^\l_n}\bigr]^{4/10}.
\end{eqnarray*}

We already know from assumption (\ref{eq:varass}) that
$\frac{1}{(\l n)^2} \esp^0[(\vert Z_n\vert-n\ell(\mu_0;d))^2]$ is
bounded in $n$.
Let us prove that $\esp^0[e^{3\l M_n-3\l^2 A_n+3 R^\l_n}]$ is also
bounded in $n$.

Note that there exists a constant $C$ such that
\begin{eqnarray*}
\esp^0\bigl[e^{3\l M_n-3\l^2 A_n+3 R^\l_n}\bigr]&\le&\esp ^0
\bigl[e^{3\l M_n+3
R^\l_n}\bigr]
\\
&=&e^{3\esp^0[R^\l_n]}\esp^0\bigl[e^{3\l M_n+3 (R^\l_n-\esp^0[R^\l
_n])}\bigr]\\
&\le&
e^C \esp^0\bigl[e^{3\l M_n+3 (R^\l_n-\esp^0[R^\l_n])}\bigr].
\end{eqnarray*}
(We used the fact that $\esp^0[R^\l_n]$ is bounded for the last inequality.)

From the independence of the $X_j$'s follows that
\[
\esp^0\bigl[e^{3\l M_n+3 (R^\l_n-\esp^0[R^\l_n])}\bigr] =\esp^0
\bigl[e^{3\l\nu(X_1)+3\l^2(Y^\l_1-\esp^0[Y^\l_1])}\bigr]^n.
\]
But the random variables $\nu(X_1)+\l(Y^\l_1-\esp^0[Y^\l_1])$ are
centered and\break bounded (uniformly in $\l$), that is, there exists a
number $M$ such that $\vert\nu(X_1)+\l(Y^\l_1-\esp^0[Y^\l
_1])\vert\le
M$ for all $\l$ and all trajectory $\omega$. Therefore,\break
$\esp^0[e^{3\l\nu(X_1)+3\l^2(Y^\l_1-\esp^0[Y^\l_1])}]^n$ is bounded
whenever $\l^2n$ is also bounded. For this last step, we rely on the
following classical lemma; see part 7 of \cite{kn:dob}, for instance.

%le3.7 #&#
\begin{lm} \label{lem:dobr}
For all $M$ and $K$, there exist constants $C_M$ and $n_0$ s.t. for all
random variable $X$
with $\vert X\vert\le M$ and $\esp[X]=0$ and for all $\l$ and $n\ge
n_0$ s.t. $\l^2n\le K$ then
\[
\esp\bigl[e^{\l X}\bigr]^n\le e^{C_M}.
\]
\end{lm}

\begin{pf}
Since $X$ is bounded, the log-Laplace transform
\[
\Lambda(\l):=\log\esp\bigl[e^{\l X}\bigr]
\]
is analytic in a neighborhood of $0$.
More precisely, we have:

Let $y\in\mathbb{C}$. Then
\[
\bigl\vert e^y-1-y\bigr\vert\le\bigl\vert e^{\vert y\vert}-1-\vert y\vert\bigr\vert\le
\vert y\vert^2e^{\vert y\vert}.
\]
Apply this to $\l X$, $\l\in\mathbb{C}$:
\[
\bigl\vert\esp\bigl[e^{\l X}\bigr]-1\bigr\vert\le\vert\l\vert^2M^2e^{\vert\l\vert
M}
\le \tfrac{1}2,
\]
if $\vert\l\vert\le\l_0$.
So $\vert\esp[e^{\l X}]\vert\ge\frac{1}2$ for $\vert\l\vert\le\l_0$.
So $\Lambda$ is analytic in $\{\l \mbox{ s.t. } \vert\l\vert\le\l_0\}$ and
$\vert\Lambda(\l)\vert\le c_0$ for some constant $c_0$, for all $\l$
s.t. $\vert\l\vert\le\l_0$. Note that
$\l_0$ and $c_0$ depend only on $M$.

We have $\Lambda(0)=0$ and $\Lambda'(0)=0$ (because $X$ is centered).
So the function $\l^{-2} \Lambda(\l)$ is also analytic in $\{\l \mbox{ s.t. }
\vert\l\vert\le\l_0\}$.
By the maximum principle, for any $\l$ such that $\vert\l\vert\le
\l_0$,
we have
\[
\bigl\vert\Lambda(\l)\bigr\vert\le C \vert\l\vert^2,
\]
where $C=\max_{z; \vert z\vert=\l_0}\frac{\Lambda(z)}{z^2}\le
\frac
{c_0}{\l_0^2}$.

The statement of the lemma is thus proved with $n_0$ chosen such that
$K/n_0\le\l_0^2$ and $C_M=c_0 K/\l_0^2$.

This completes the proof of (\ref{eq:fond'}) in the case $\alpha\neq0$.

The case $\alpha=0$ is easier. As we did in (\ref{eq:fond0}) and
(\ref
{eq:fond1}), we start writing that
\[
\frac{1}{\l n} \bigl(\esp^\l\bigl[\vert Z_n\vert\bigr]-
\esp^0\bigl[\vert Z_n\vert\bigr]\bigr) =\frac{1}{\l n}
\esp^0\bigl[\bigl(\vert Z_n\vert-n\ell(
\mu_0;d)\bigr) \bigl(e^{\l M_n-\l^2
A_n+ R^\l_n}-1\bigr)\bigr].
\]
Using similar arguments as for the case $\alpha\neq0$, it is not
difficult to show that
\begin{eqnarray*}
&&\lim_{n\rightarrow+\infty,\l\rightarrow0,\l^2n\ra0}\frac{1}{\l n} \esp^0\bigl[\bigl(
\vert Z_n\vert-n\ell(\mu_0;d)\bigr) \bigl(e^{\l M_n-\l^2 A_n+ R^\l
_n}-1
\bigr)\bigr]
\\
&&\qquad= \lim_{n\rightarrow+\infty,\l\rightarrow0,\l^2n\ra0}\frac{1}{\l n} \esp^0\bigl[\bigl(
\vert Z_n\vert-n\ell(\mu_0;d)\bigr) \l M_n
e^{-\l^2 A_n+ R^\l
_n}\bigr].
\end{eqnarray*}

%Use the bound $\vert e^y-1-y\vert\le y^2(1+e^y)$. Then use Holder with
%exponent $2$ on $\vert Z_N\vert$, exponent $6$ for $(\l M_n)^2$,
%exponent $6$ for $(1+e^{\l M_n})$ and exponent $6$ for $e^{ R^\l_n}$.
%(The term $A_n$ disappears because of its sign.) We get four terms of
%respective orders $\sqrt n$, $\l^2n$, $1$ and $1$ (using estimates on
%the log-Laplace of $\l M_n$
%and $R^\l_n$ as before). Note that $\sqrt{n}\l^2n<<\l n$.
%

Observe that
\begin{eqnarray*}
&&\frac{1}{\l n} \esp^0\bigl[\bigl(\vert Z_n\vert-n
\ell(\mu_0;d)\bigr) \l M_n e^{-\l^2
A_n+ R^\l_n}\bigr] \\
&&\qquad=
\frac{1}{n} \esp^0\bigl[\bigl(\vert Z_n\vert-n\ell(
\mu_0;d)\bigr) M_n e^{-\l^2 A_n+
R^\l_n}\bigr].
\end{eqnarray*}

The limit of this last expression is given by the central limit theorem
and, with the notation we already used, it coincides with $\esp
[ZM]=\sigma$. Observe that,
with our scaling satisfying $\l^2n\ra0$, we have $\l^2 A_n\ra0$.
The details are similar to the case $\alpha\neq0$.
\end{pf}

\begin{pf*}{End of the proof of Theorem~\ref{theo:rescp}}
Proposition~\ref{prop:clt} and Lemma~\ref{lem:2ndmoment} from Section
\ref
{sec:centrallimittheorems} show that the assumptions of Theorem~\ref
{theo:malliavin} are satisfied. Thus, we get that
\[
\lim_{n\rightarrow+\infty,\l\ra0, \l^2n\ra1} \frac{1}{\l n} \bigl(\esp ^\l \bigl[\vert
Z_n\vert\bigr]-\esp^0\bigl[\vert Z_n\vert\bigr]\bigr) =
\sigma(\nu,\mu_0;d).
\]
But we observed in Section \ref{sec:geometricinput} that this convergence
implies Theorem~\ref{theo:rescp}.
\end{pf*}

%%%%%%%%%%%%%%%%%%%%%%%%%%%%%%%%%%%%%%%%%%%%%%%%%%%%%%%%%%%%%%%%%%%%%%%%%%%%%%%%%%%%%
%%%%%%%%%%%%%%%%%%%%%%%%%%%%% ENTROPY
%%%%%%%%%%%%%%%%%%%%%%%%%%%%%%%%%%%%%%%%%%%%%%%
%%%%%%%%%%%%%%%%%%%%%%%%%%%%%%%%%%%%%%%%%%%%%%%%%%%%%%%%%%%%%%%%%%%%%%%%%%%%%%%%%%%%%

%s4 #&#
\section{Proof of Theorem \texorpdfstring{\protect\ref{theo:entrop}}{2.2}}
\label{sec:entrop}

We now explain how to deduce Theorem~\ref{theo:entrop} from Theorem~\ref
{theo:rescp}.

As in the proof of Theorem~\ref{theo:rescp}, we may and will restrict
ourselves to positive~$\l$'s.

We first recall that the entropy can be interpreted as a rate of escape
in the Green metric: $h(\mu)=\ell(\mu;d_G^\mu)$.
Thus, we have
%
%e4.1 #&#
\begin{equation}
\label{eq:entrop1} h(\mu_\l)-h(\mu_0)= \bigl(h(
\mu_\l)-\ell\bigl(\mu_\l;d_G^0
\bigr) \bigr) + \bigl(\ell\bigl(\mu _\l;d_G^0
\bigr)-\ell\bigl(\mu_0;d_G^0\bigr) \bigr).
\end{equation}
By Theorem~\ref{theo:rescp} and since $d_G^0$ satisfies (\emph{BA}), once
divided by $\l$, the second term in~(\ref{eq:entrop1}) converges to
$\sigma(\nu,\mu_0;d_G^0)$. Thus, the proof of Theorem~\ref{theo:entrop}
will be complete once we prove that
%
%e4.2 #&#
\begin{equation}
\label{eq:entroptrop} \lim_{\l\ra0+} \frac{1} \l \bigl(h(
\mu_\l)-\ell\bigl(\mu_\l;d_G^0
\bigr) \bigr)=0.
\end{equation}

It turns out the estimate (\ref{eq:entroptrop}) does not use the
hyperbolicity of $\G$. We have the following more general property.

%pr4.1 #&#
\begin{prop} \label{prop:prop}
Let $\G$ be a finitely generated group. Assume $\G$ is nonamenable.
Let $\mu_0$ be a probability measure on $\G$ such that
the support of $\mu_0$ generates $\G$ (as a semigroup), and
$H(\mu_0)<\infty$.

Consider a curve of probability measures on $\G$, say $\l\in
[-1,1]\ra\mu
_\l$, satisfying the Regularity Assumption:
\[
\log\mu_\l(a)=\log\mu_0(a)+\l\nu(a)+\l
o_\l(a),
\]
where
$\nu$ is bounded and
$o_\l(a)$ converges to $0$ uniformly in $a\in\G$.

Then
\[
\lim_{\l\ra0+} \frac{1} \l \bigl(h(\mu_\l)-\ell
\bigl(\mu_\l;d_G^0\bigr) \bigr)=0.
\]
\end{prop}

%re4.2 #&#
\begin{rmk}
We do not assume any more that $\mu_0$ or the $\mu_\l$'s are symmetric.
Then the Green metric may not be a real distance. Indeed, although it
still satisfies the triangle inequality, it may not be symmetric.

Thus, in this part of the paper, the word ``metric'' will refer to a
function on $\G\times\Gamma$ that vanishes on the diagonal and
satisfies the triangle inequality.

The interpretation of the asymptotic entropy as the rate of escape in
the Green metric remains valid in this general framework; see \cite{kn:bhm1}.
\end{rmk}
%

%re4.3 #&#
\begin{rmk}
The assumption that $H(\mu_0)<\infty$ implies that $\mu_0$ has a finite
first moment with respect to $d_G^0$
(see Lemma~2.3 in \cite{kn:bhm1}).
The Regularity Assumption then implies that $H(\mu_\l)<\infty$ and that
$\mu_\l$ also has a finite first moment with respect to $d_G^0$.
\end{rmk}

%Proof shows that $$\vert h(\mu')-\ell(\mu';d_G)\vert\leq\frac c{
%provide $\mu'(a)/\mu(a)$ is uniformly small. \end{rmk}{\chr"CA}

\begin{pf*}{Proof of Proposition \ref{prop:prop}}
We give two separate arguments for lower and upper bounds for $h(\mu
_\l
)-\ell(\mu_\l;d_G^0)$. Both arguments are based on the so-called
``fundamental inequality'' that we first recall:
let $\mu$ be a probability measure on $\G$ with finite entropy and let
$d$ a left-invariant proper metric on $\G$. We denote with ${\mbox v}(d)$
the logarithmic volume growth of the metric $d$.
The ``fundamental inequality'' states that $h(\mu)\leq{\mbox v}(d) \ell
(\mu;d)$; see \cite{kn:guiv,kn:versh,kn:bhm1} and the
references quoted therein.

The ``fundamental inequality'' in particular applies to any Green
metrics $d_G^\a$. By a result in \cite{kn:bhm1}, we have ${\mbox
v}(d_G^\a)=1$. Therefore, we get that
%
%e4.3 #&#
\begin{equation}
\label{eq:ineqfund} h(\mu_\l)=\ell\bigl(\mu_\l;d_G^\l
\bigr)\leq\ell\bigl(\mu_\l;d_G^\a\bigr)
\end{equation}
for all $\l$ and $\a$.

Applying (\ref{eq:ineqfund}) with $\alpha=0$, yields $h(\mu_\l)\leq
\ell(\mu_\l;d_G^0)$ and, therefore,
%
%e4.4 #&#
\begin{equation}
\label{eq:entropup} \limsup_{\l\ra0+} \frac{1} \l \bigl(h(
\mu_\l)-\ell\bigl(\mu_\l;d_G^0
\bigr) \bigr)\leq 0.
\end{equation}

It remains to prove the lower bound. We first need review properties of
the Green metric.
Consider a probability measure $\mu$ with finite entropy and whose
support generates the whole group $\Gamma$.
We recall that we defined the Green metric as
\[
d_G^\mu(x,y):=\log G^\mu(\mathrm{id})-\log
G^\mu\bigl(x^{-1}y\bigr),
\]
where $G^\mu(x)=\sum_{n=0}^\infty\mu^n(x) $ is the Green function of
the random walk.

We may equivalently express $d_G^\mu$ in terms of the hitting
probabilities of the random walk:
for a given trajectory $\omega\in\Omega$ and $z\in\G$, let
\[
T_z(\omega)=\inf\bigl\{n\geq0; Z_n(\omega)=z\bigr\}
\]
be the hitting time of $z$ by $\omega$. Observe that $T_z(\omega)$
may be infinite.

Define $F^\mu(z):=\P^\mu[T_z<\infty]$. Then
\[
d_G^\mu(\mathrm{id},z)=-\log F^\mu(z)
\]
as can be easily checked using the Markov property.

In the sequel, we use the notation $F^\l$ instead of $F^{\mu_\l}$.

%le4.4 #&#
\begin{lm}\label{lm:ub}
The function $(\l,\a)\ra\ell(\mu_\l;d_G^\a)$ is bounded on $[0,1]^2$:
\[
\sup_{0\leq\l\leq1; 0\leq\a\leq1}\ell\bigl(\mu_\l;d_G^\a
\bigr)<\infty.
\]
\end{lm}

 \begin{pf}
Let $\mu$ and $\mu'$ be probability measures on $\G$ and let $d$ be a
proper left-invariant metric on $\G$.
It is clear that
\[
\ell(\mu;d)\leq\sum_{a\in\G} d(\mathrm{id},a) \mu(a).
\]

Also $F^{\mu'}(a)\geq\mu'(a)$ and, therefore, $d_G^{\mu'}(\mathrm{id},a)\leq
-\log\mu'(a)$.

Applying these two inequalities to $\mu=\mu_\l$ and $\mu'=\mu
_\alpha$,
we get that
\begin{eqnarray*}
\ell\bigl(\mu_\l;d_G^\alpha\bigr)&\leq&-\sum
_{a\in\G}\bigl(\log\mu_\alpha (a)\bigr)
\mu_\l (a)
\\
&=&-\sum_{a\in\G}\bigl(\log\mu_0(a)+\alpha
\nu(a)+\alpha o_\alpha (a)\bigr)\mu_\l (a).
\end{eqnarray*}
Since $\nu$ and $o_\l$ are bounded, we have $\mu_\l(a)\leq e^C\mu_0(a)$
for some constant $C$.
For the same reason, the term $\sum_{a\in\G}(\alpha\nu(a)+\alpha
o_\alpha(a))\mu_\l(a)$ is also controlled by a constant.
Therefore,
\[
\ell\bigl(\mu_\l;d_G^\alpha\bigr)\leq
e^CH(\mu_0)+C
\]
for some constant $C$.
\end{pf}

We shall need the following estimate on $T_z$.

%le4.5 #&#
\begin{lm}\label{lm:tau}
Let $\mu$ be a probability measure on $\Gamma$ whose support
generates~$\Gamma$.
Then there exists a positive constant
$\kappa$
such that
\[
\sup_{z\in\G} \esp^\mu\bigl[e^{\kappa T_z};
T_z<\infty\bigr]<\infty.
\]
\end{lm}

\begin{pf}
We use the nonamenability of $\G$: there exists a constant $\rho<1$
such that, for all $n$ and all $z\in\G$, we have
$\mu^n(z)\leq\rho^n$. Therefore,
\begin{eqnarray*}
\esp^\mu\bigl[e^{\kappa T_z}; T_z<\infty\bigr] &=&
\sum_n e^{\kappa n} \P^\mu[T_z=n]
\\
&\leq& \sum_n e^{\kappa n}
\P^\mu[Z_n=z]=\sum_n
e^{\kappa n} \mu ^n(z)
\\
&\leq& \sum_n e^{\kappa n}
\rho^n<\infty
\end{eqnarray*}
as soon as $e^\kappa\rho<1$.
\end{pf}

\begin{pf*}{End of the proof of Proposition~\ref{prop:prop}: The lower bound}
We use the shorthand notation $F^\l(z):=F^{\mu_\l}(z)$.

As in (\ref{eq:fond0}), we have
\[
F^\l(z)=\esp^0\Biggl[\prod_{j=1}^{T_z}
\frac{\mu_\l(X_j)}{\mu_0(X_j)}; T_z<\infty\Biggr].
\]
Let us apply H\"older's inequality with positive parameters $(p,q,r)$
such that $1/p+1/q+1/r=1$ and with the notation
$\a=\l q$. We assume that $\a\leq1$. Thus,
%
%e4.5 #&#
\begin{eqnarray}
\label{eq:f1} F^\l(z)&\leq& F^0(z)^{1/p}
F^\a(z)^{1/q}
\nonumber
\\[-8pt]
\\[-8pt]
\nonumber
&&{}\times \esp^0\Biggl[ \Biggl(
\prod_{j=1}^{T_z} \frac{\mu_\l(X_j)}{\mu
_0(X_j)}\biggl(\frac{\mu
_0(X_j)}{\mu_\a(X_j)}
\biggr)^{1/q} \Biggr)^r; T_z<\infty
\Biggr]^{1/r}.
\end{eqnarray}

Let $a\in S$.
Using the Regularity Assumption and the relation $\a=\l q$, we get
\[
\frac{\mu_\l(a)}{\mu_0(a)}\biggl(\frac{\mu_0(a)}{\mu_\a(a)}\biggr)^{1/q}
=e^{\l(o_\l(a)-o_\a(a))}.
\]
Therefore, since $\l\leq\a\leq1$,
and remembering that $o_\l$ uniformly converges to $0$, we see that for
all $\eps>0$, provided $\a$ is small enough then
%
%e4.6 #&#
\begin{equation}
\label{eq:f2} \frac{\mu_\l(a)}{\mu_0(a)}\biggl(\frac
{\mu
_0(a)}{\mu_\a(a)}\biggr)^{1/q}
\leq\exp(\eps\l).
\end{equation}

Using (\ref{eq:f2}) in equation (\ref{eq:f1}), we see that
\[
F^\l(z)\leq F^0(z)^{1/p} F^\a(z)^{1/q}
\esp^0\bigl[e^{\eps\l rT_z}; T_z<\infty\bigr].
\]
If we further assume that $\eps\l r\leq\kappa^0$, where $\kappa^0$ is
the constant given by Lemma~\ref{lm:tau} when choosing $\mu=\mu_0$,
then we have
%
%e4.7 #&#
\begin{equation}
\label{eq:f3} F^\l(z)\leq C F^0(z)^{1/p}
F^\a(z)^{1/q}
\end{equation}
for a new constant $C$ that does not depend on $z$.

We evaluate inequality (\ref{eq:f3}) at $z=Z_n$; take the logarithm and
take the expectation with respect to $\P^\l$ to obtain
\[
\esp^\l\bigl[d_G^\l(\mathrm{id},Z_n)
\bigr]\geq\frac{1} p \esp^\l\bigl[d_G^0(\mathrm{id},Z_n)
\bigr]+\frac{1} q \esp^\l\bigl[d_G^\a(\mathrm{id},Z_n)
\bigr]-\log C.
\]
Now divide by $n$ and let $n$ tend to $\infty$, to get that
%
%e4.8 #&#
\begin{equation}
\label{eq:f4} h(\mu_\l)=\ell\bigl(\mu_\l;d_G^\l
\bigr)\geq\frac{1} p \ell\bigl(\mu_\l;d_G^0
\bigr)+\frac{1} q \ell\bigl(\mu_\l;d_G^\a
\bigr).
\end{equation}
We choose $r=\kappa^0/(\eps\l)$ and $\l$ small enough so that
$1/r+1/q<1$. Then (\ref{eq:f4}) becomes
%
%e4.9 #&#
\begin{equation}
\label{eq:f5} h(\mu_\l)-\ell\bigl(\mu_\l;d_G^0
\bigr) \geq\frac\l\a \bigl(\ell\bigl(\mu_\l;d_G^\a
\bigr)-\ell\bigl(\mu_\l;d_G^0\bigr) \bigr)-
\frac
{\eps\l}{\kappa^0} \ell\bigl(\mu_\l;d_G^0
\bigr).
\end{equation}

We let $\l$ tend to $0$ in (\ref{eq:f5}): by Lemma~\ref{lm:ub}, we
know that
\[
\l \bigl(\ell\bigl(\mu_\l;d_G^\a\bigr)-\ell
\bigl(\mu_\l;d_G^0\bigr) \bigr)\ra0
\]
and $\l\ell(\mu_\l;d_G^0)\ra0$. Therefore,
%
%e4.10 #&#
\begin{equation}
\label{eq:f55} \liminf_{\l\ra0+} \bigl(h(\mu_\l)-\ell
\bigl(\mu_\l;d_G^0\bigr) \bigr)\geq 0.
\end{equation}

Using the inequality $\ell(\mu_\l;d_G^\a)\geq h(\mu_\l)$ (which comes
from the ``fundamental inequality''), we deduce from (\ref{eq:f5}) that
%
%e4.11 #&#
\begin{equation}
\label{eq:f56} \frac{1}\l{ \bigl(h(\mu_\l)-\ell\bigl(
\mu_\l;d_G^0\bigr) \bigr)} \geq
\frac{1}\a \bigl(h(\mu_\l)-\ell\bigl(\mu_\l;d_G^0
\bigr) \bigr)-\frac
{\eps}{\kappa
^0} \ell\bigl(\mu_\l;d_G^0
\bigr).
\end{equation}
It follows from Lemma~\ref{lm:ub} that there exists a constant $\ell_0$
such that\break
$\frac{1}{\kappa^0} \ell(\mu_\l;d_G^0)\leq\ell_0$ for all $\l$.
By (\ref{eq:f55}), the term $h(\mu_\l)-\ell(\mu_\l;d_G^0)$ has a
nonnegative $\liminf$. Thus, we deduce from (\ref{eq:f56}) that
%-c\frac{\a}{\kappa^0} \ell(\mu_0;d_G^0).\eeqn
%By the 'fundamental inequality', we know that $\ell(\mu_0;d_G^\a)\geq
%
%e4.12 #&#
\begin{equation}
\label{eq:f7} \liminf_{\l\ra0+} \frac{1}\l \bigl(h(
\mu_\l)-\ell\bigl(\mu_\l;d_G^0
\bigr) \bigr)\geq -\eps\ell_0.
\end{equation}
And since (\ref{eq:f7}) holds for any small enough $\eps$, we have
%
%e4.13 #&#
\begin{equation}
\label{eq:entropdown} \liminf_{\l\ra0+} \frac{1}\l \bigl(h(
\mu_\l)-\ell\bigl(\mu_\l;d_G^0
\bigr) \bigr)\geq 0.
\end{equation}
\upqed\end{pf*}
\noqed\end{pf*}

%re4.6 #&#
\begin{rmk}
F. Ledrappier and L. Shu recently adapted our strategy in a continuous setting:
using martingales as here, they obtained differentiability results for
the entropy and rate of escape
of Brownian motions on the universal cover of negatively curved manifolds;
see \url{http://front.math.ucdavis.edu/1309.5182}.
\end{rmk}

%%%%%%%%%%%%%%%%%%%%%%%%%%%%%%%%%%%%%%%%%%%%%%%%%%%%%%%%%%%%%%%%%%%%%%%%%%%%%%%%%%%%%
%%%%%%%%%%%%%%%%%%%%%%%%%%%%% MISC
%%%%%%%%%%%%%%%%%%%%%%%%%%%%%%%%%%%%%%%%%%%%%%%
%%%%%%%%%%%%%%%%%%%%%%%%%%%%%%%%%%%%%%%%%%%%%%%%%%%%%%%%%%%%%%%%%%%%%%%%%%%%%%%%%%%%%

%
%
%For any Gaussian vector $(Z,M)$
%$$\esp[Ze^{M-\frac1 2\esp[M^2]}]
%=\esp[ZM].$$
%
%Let $\Sigma$ be the covariance matrix of the vector $(Z,M)$.
%
%For $x=(\l,\a)\in\R^2$, let $(x,x)_\Sigma$ be the square of the norm
%of $x$ in the metric induced by $\Sigma$.
%We have
%$$\esp[e^{\l Z+\a M}]=e^{\frac12 (x,x)_\Sigma}.$$
%Choose $\a=1$ and take the derivative at $\l=0$.\qed
%
%If $(\rho_n)$ weakly converges to $\rho$ and $\sup_n\int\phi^{1+\a} d
%$$\int\phi d\rho_n\ra\int\phi d\rho.$$
%
%
%

%%%%%%%%%%%%%%%%%%%%%%%%%%%%%%%%%%%%%%%%%%%%%%%%%%%%%%%%%%%%%%%%%%%%%%%%%%%%%%%%%%%%%
%%%%%%%%%%%%%%%%%%%%%%%%%%%%% BIB
%%%%%%%%%%%%%%%%%%%%%%%%%%%%%%%%%%%%%%%%%%%%%%%%%
%%%%%%%%%%%%%%%%%%%%%%%%%%%%%%%%%%%%%%%%%%%%%%%%%%%%%%%%%%%%%%%%%%%%%%%%%%%%%%%%%%%%%

% imsref loaded by akundreckaite, 2014-02-06 07:59:37
% imsref loaded by akundreckaite, 2014-02-06 08:00:56
%

%


%suskaldyti doi

\printaddresses

\end{document}